\documentclass[12pt]{article}
\usepackage{amsfonts} \input amssym.def \input amssym

\def\FF{{\mathbb F}}  \def\NN{{\mathbb N}}
  
\newtheorem{theorem}{Theorem} 
\newtheorem{pr}{Property} \newtheorem{defi}{Definition}
 \newtheorem{lemma}{Lemma}

\def\endproof{\nobreak\kern5pt\nobreak\vrule height4pt width4pt
depth0pt \vskip4pt plus2pt}

\newenvironment{matrice}{\left(\begin{array}}{\end{array}\right)}

\begin{document}
\title{Skew-cyclic codes} 

\author{D.~Boucher\thanks{IRMAR, Universit\'e de Rennes 1,
Campus de Beaulieu, F-35042 Rennes Cedex}, W.~Geiselmann\thanks{IAKS, Universit{\"a}t Karlsruhe, Fakult{\"a}t f\"ur Informatik, Postfach 6980, D-76128 Karlsruhe} and  F.~Ulmer \thanks{IRMAR, Universit\'e de Rennes 1,
Campus de Beaulieu, F-35042 Rennes Cedex}}

\date{April 17, 2006}

\maketitle

\begin{abstract} We generalize the notion of cyclic codes by using generator polynomials in (non commutative) skew polynomial rings. Since skew polynomial rings are left and right euclidean, the obtained codes share most properties of cyclic codes. Since there are much more skew-cyclic codes, this new class of codes allows to systematically search for codes with good properties. We give many examples of codes which improve the previously best known linear codes.
\end{abstract} 

\section*{Introduction}
Let $\FF_q$ be a finite field of $q$ elements.  A linear $(n,k)$-code over $\FF_q$ is a $k$-dimensional vector subspace ${\cal C}$ of the vector
space $$V={\FF_q}^n=\{(a_0, \ldots, a_{n-1})\, | \, a_i \in \FF_q\}.$$
In the following we use the polynomial representation of the code.  In
the polynomial representation of the code ${\cal C}$,
the code words $(a_0,a_1,\ldots,a_{n-1})\in {\cal C}$ are coefficient
tuples of elements
$a_{n-1}{X}^{n-1}+\ldots+a_1{X}+a_0\in\FF_q[X]/(X^n-1)$ which are
multiples of one element $G\in A$ (the generator polynomial).  A
linear code  ${\cal C}$ is a cyclic code if \[
(a_0,a_1,\ldots,a_{n-1})\in {\cal C} \quad \Rightarrow \quad
(a_{n-1},a_0,a_1,\ldots,a_{n-2})\in {\cal C}.  \]

In this paper we want to generalize the notion of cyclic codes to the
notion of $\theta$-cyclic codes.

\begin{defi} Let $\FF_q$ be a finite field and $\theta$ an
automorphism of $\FF_q$.  A $\theta$-cyclic code is a linear code
${\cal C}_{\theta}$ with the property that \[
(a_0,a_1,\ldots,a_{n-1})\in {\cal C}_{\theta} \quad \Rightarrow \quad
(\theta(a_{n-1}),\theta(a_0),\theta(a_1),\ldots,\theta(a_{n-2}))\in
{\cal C}_{\theta}.  \]
\end{defi}

In order to generalize the notion of cyclic codes (corresponding to
the case where $\theta$ is the identity) we consider skew polynomial
rings of automorphism type which we now define.  Starting from the
finite field $\FF_q$ and an automorphism $\theta$ of $\FF_q$ one
defines a ring structure on the set \[ \FF_q[X,\theta]\, = \, \left\{
a_{n}{X}^{n-1}+\ldots+a_1{X}+a_0 \, | \, a_i\in \FF_q \mbox{ and }
n\in \NN\right\}.  \] This is the set of formal polynomials where the
coefficients are written on the left of the variable $X$.  The
addition in $\FF_q[X,\theta]$ is defined to be the usual addition
of polynomials and the multiplication is defined by the basic rule
$Xa=\theta(a)X$ ($a\in \FF_q$) and extended to all elements of
$\FF_q[X,\theta]$ by associativity and distributivity.  Those rings
are well known (cf.\ \cite{OR,Mc}) and, since over a finite field all
derivations are inner, they are the most general ``polynomial
rings'' with a commutative field of coefficients where the degree of a
product of two elements is the sum of the degrees of the elements.

Our goal is to give a skew polynomial representation of
$\theta$-cyclic codes.  We will show that the code words
$(a_0,a_1,\ldots,a_{n-1})$ of a $\theta$-cyclic code ${\cal
C}_{\theta}$ are coefficient tuples of elements of
$a_{n-1}{X}^{n-1}+\ldots+a_1{X}+a_0\in \FF_q[X,\theta]/(X^n-1)$ which
are left multiples of one element $G\in \FF_q[X,\theta]/(X^n-1)$
(the generator polynomial).  This property also guaranties that the
encoding procedure of a $\theta$-cyclic code is as easy as for cyclic
codes.

 We will also show by concrete examples that the class of
 $\theta$-cyclic codes is a very large  class of linear codes (containing the cyclic 
 codes) and that this class contains codes with good
 properties.  Therefore the class $\theta$-cyclic codes is an interesting class
 of linear codes which are easy to construct in a systematic way.  In
 a final section we will show how to decode some $\theta$-cyclic codes.

\section{Generalities on $\theta$-cyclic codes}

Properties of $\theta$-cyclic codes are closely related to properties
of $\FF_q[X,\theta]$.  The ring $\FF_q[X,\theta]$ is a left and right
euclidean ring whose left and right ideals are principal \cite{OR}.
Here right division means that for $P_1(X),P_2(X)\in \FF_q[X,\theta]$
which are non zero, there exist unique polynomials $Q_r(X),R_r(X)\in
\FF_q[X,\theta]$ such that \[ P_1(X)=Q_r(X)\cdot P_2(X)\, +\,
P_r(X).\] If $P_r(X)=0$ then $P_2(X)$ is a right divisor of $P_1(X)$
in $\FF_q[X,\theta]$.  The definition of left divisor in
$\FF_q[X,\theta]$ is similar using the left euclidean division.  In
the ring $\FF_q[X,\theta]$ left and right gcd and lcm exist and can be
computed using the left and right euclidean algorithm.  We denote
${\cal F}\subset \FF_q$ the subfield of elements of $\FF_q$ that are
left fixed by $\theta$.  An element $P\in \FF_q[X,\theta]$ is central
(i.e.\ commutes with all elements of $\FF_q[X,\theta]$) if and only if
$P=\sum_{i=0}^m c_iX^{i\cdot \alpha}\in {\cal F}[X]$ where
$\alpha=|<\theta>|$ is the order of $\theta$ (\cite{Mc}, Theorem
II.12).  In particular central elements of $\FF_q[X,\theta]$ are the
generators of two-sided ideals in $\FF_q[X,\theta]$. Therefore, if
$|<\theta>|$ divides $n$, then $(X^n-1)\subset \FF_q[X,\theta]$ is a
two-sided ideal.  In the non-commutative ring
$\FF_q[X,\theta]/(X^n-1)$ we identify the image of $P\in
\FF_q[X,\theta]$ under the canonical morphism $ \psi\colon
\FF_q[X,\theta] \to \FF_q[X,\theta]/(X^n-1)$ with the remainder $P_r$
of $P$ by the right division with $X^n-1$ in $\FF_q[X,\theta]$ and we
denote $\psi(X)$ still by $X$.  This representation gives a canonical
form for the elements of $\FF_q[X,\theta]/(X^n-1)$.

\begin{lemma} Let $\FF_q$ be a finite field, $\theta$ an automorphism
of $\FF_q$ and $n$ an integer divisible by the order $|<\theta>|$ of
$\theta$.  The ring $\FF_q[X,\theta]/(X^n-1)$ is a principal left
ideal domain in which left ideals are generated by $\psi(G)$ where $G$
is a right divisor of $X^n-1$ in $\FF_q[X,\theta]$.
\end{lemma}
\begin{proof} The proof is an exact copy of the commutative case only
taking care of left and right.  Let $I$ be a left ideal of
$\FF_q[X,\theta]/(X^n-1)$.  If $I=\{0\}$ then $I=(0)$.  Otherwise
denote $G_r\in I$ a monic non zero polynomial of minimal degree in
$I$.  By abuse of notation we identify the element $G_r\in I$ with
itself in $\FF_q[X,\theta]$ and denote this element $G$ (i.e.~$G\in
\FF_q[X,\theta]$ is of degree $<n$ and $\psi(G)=G_r$).  Let $P_r\in I$
be an arbitrary element of $I$ which we again identify with $P\in
\FF_q[X,\theta]$.  Performing a right division of $P$ by $G$ in
$\FF_q[X,\theta]$ we get \[P\, =\, Q\cdot G \, + R,\qquad \mbox{ where
} \mbox{deg}(R)<\mbox{deg}(G)\] from which we get
$\psi(R)=P_r-\psi(Q)\cdot G_r \in I$.  By minimality of the degree of
$G_r$ we must have $\psi(R)=0$, showing that $P_r=\psi(Q)\cdot G_r$
and thus $I=(G_r)$.
\end{proof}

For a linear code ${\cal C}$ of length $n$ we denote ${\cal C}(X)$ the
skew polynomial representation of ${\cal C}$.  In this representation
we associate to a code word $a=(a_0,a_1,\ldots,a_{n-1})\in {\cal C}$
the element $a(X)=a_{n-1}{X}^{n-1}+\ldots+a_1{X}+a_0$ in
$\FF_q[X,\theta]/(X^n-1)$.  If $a\in {\cal C}$, then we denote
$a(X)\in \FF_q[X,\theta]/(X^n-1)$ the skew polynomial representation
of $a$.

\begin{theorem} Let $\FF_q$ be a finite field, $\theta$ an
automorphism of $\FF_q$ and ${\cal C}$ be a linear code over $\FF_q$
of length $n$.  If $|<\theta>|$, the order of $\theta$,  divides $n$, then the code ${\cal C}$
is a $\theta$-cyclic code if and only if the skew polynomial
representation ${\cal C}(X)$ of ${\cal C}$ is a left ideal $(G)
\subset \FF_q[X,\theta]/(X^n-1)$.
\end{theorem}
\begin{proof} By the above Lemma we have to show that ${\cal C}(X)$ is a
left ideal of $\FF_q[X,\theta]/(X^n-1)$.  Since ${\cal C}$ is a linear
code, ${\cal C}(X)$ is an additive group.  Let
$a=(a_0,\ldots,a_{n-1})\in {\cal C}$, then
\begin{eqnarray*}
X \, a(X) & = & X \, a_0 + X \, (a_1 \, X) + \cdots + X \, (a_{n-1} \,
X^{n-1})\\
 & = & \theta(a_0) \, X + (\theta(a_1) \, X) \, X + \cdots +
 (\theta(a_{n-1}) \, X) \, X^{n-1}\\
 & = & \theta(a_{n-1})+\theta(a_0) \, X + \cdots \theta(a_{n-2}) \,
 X^{n-1} + \theta(a_{n-1}) \, (X^n-1).
\end{eqnarray*}
Therefore in $\FF_q[X,\theta]/(X^n-1)$ (i.e.~working modulo $X^n-1$)
we have $X \, a(X)=\theta(a_{n-1})+\theta(a_0) \, X + \cdots
\theta(a_{n-2}) \, X^{n-1}$.  Since ${\cal C}$ is $\theta$-cyclic we have $X \,
a(X)\in {\cal C}(X)$ and by iteration and linearity we get  for all
$P_r\in \FF_q[X,\theta]/(X^n-1)$ that $P_r\cdot a(X) \in {\cal C}(X)$.
This shows that ${\cal C}(X)$ is a left ideal of
$\FF_q[X,\theta]/(X^n-1)$.\\
In  the opposite direction the properties of a left ideal show  that the coefficient vectors of the elements of a   left ideal  $(G)
\subset \FF_q[X,\theta]/(X^n-1)$ form a linear subspace and from   $a(X)\in (G) \Rightarrow X \,
a(X)\in (G)$ we get from the above computation that the corresponding  linear code is  $\theta$-cyclic.
\end{proof}

\section{Finding good codes}

  An obvious technique for finding good linear codes (codes with a large minimum
distance $d$) is a random search.  With this technique, the probability to
find a code with better parameters than the best known codes, e.g.
according to Brouwer's table~\cite{codetable}
(http://www.win.tue.nl/~aeb/), is very small.
Many of the best known codes have some additional structure (e.g. are
cyclic codes or are constructed using cyclic codes).  Therefore a search
within the $\theta$-cyclic codes seems more promising than a random
search --- especially as, since $\FF_q[X,\theta]$ is not a unique factorization ring, there are many $\theta$-cyclic codes for a
given set of parameters $(n,k)$.

We implemented a factorization procedure in the CA-System
MAGMA\cite{magma}. This procedure 
outputs all right skew-factors of  $X^n-1$, producing the possible generator skew polynomials for  $\theta$-cyclic codes.   A right factor of degree $n-k$ of  $X^n-1$ generates a linear code with parameters $(n,k)$. If  $\theta$ is not the identity (corresponding to the cyclic codes),  then  $\FF_q[X,\theta]$ is in general not a unique factorization ring. In this case there are typically much more right factors than in the commutative  case, producing  many $\theta$-cyclic codes. Once the code is given, its minimum distance can be
calculated using the  existing MAGMA procedures.  This latter operation is very time consuming for larger codes, hence we restricted our search to smaller
codes with ground fields $\FF_{4})$ and $\FF_{9}$ and to 5000 codes in
the cases, where more skew factors for a given parameter set $(n,k)$
have been found.
With this technique we obtained a minimum distance one larger than the previously known
best code (according to Brouwer's table) for 8 parameter sets over
$\FF_{4}$. Those codes have been added to the MAGMA list of known codes.
In most cases we found many different codes with the same minimum
distance; in Table~\ref{codes4} the code parameters, the number of
codes found with these parameters (No), and a generating
polynomial for one code in this class  of parameters are given. 
\begin{table}
\def\f1{\alpha}
$$\begin{array}{r|r|p{10cm}}
    (n,k,d_{min}) & No &  \multicolumn{1}{c}{g}\\ \hline\hline
(56, 30, 14) & 1& 
$x^{26} + x^{23} + \f1\, x^{22} + \f1^2\, x^{21} + \f1\,
x^{20} + \f1^2\, x^{19} + \f1^2\, x^{18} + \f1\, x^{17} +
x^{16} + x^{14} + x^{13} + \f1\, x^{11} + \f1^2\, x^{10} +
\f1^2\, x^{9} + \f1^2\, x^{8} + \f1\, x^{7} + \f1^2\,
x^{6} + \f1\, x^{5} + \f1^2\, x^{4} + x^{2} + \f1^2\, x +
\f1^2
$\\\hline
(48,19,17) & 2&
$x^{29} + \f1^2\,
x^{28} + x^{26} + \f1\, x^{25} + \f1^2\, x^{24} + \f1\,
x^{23} + \f1\, x^{21} + \f1\, x^{20} + \f1^2\, x^{19} +
\f1\, x^{18} + \f1\, x^{17} + \f1\, x^{16} + x^{15} + x^{14}
+ \f1\, x^{13} + \f1\, x^{10} + \f1\, x^{8} + \f1^2\,
x^{7} + x^{6} + x^{5} + x^{4} + \f1^2\, x^{3} + x^{2} + \f1^2
$\\\hline
(48,25,13)& 2& 
$x^{23} + \f1^2\, x^{22} + x^{21} + \f1\, x^{20} + \f1\,
x^{19} + \f1^2\, x^{18} +\f1\, x^{17} + \f1\, x^{14} +
\f1^2\, x^{13} + \f1^2\, x^{11} + x^{9} + \f1\, x^{7} + x^{6}
+ x^{3} + \f1^2\, x^{2} + 1
$\\\hline
(42, 17, 16)& 3& 
$x^{25} + x^{23} + \f1\, x^{22} + x^{21} + x^{20} + x^{19} +
x^{18} + \f1^2\, x^{17} + \f1^2\, x^{16} + \f1\, x^{15} +
\f1\, x^{14} + x^{13} + x^{11} + x^{10} + x^{8} + \f1^2\, x^{4}
+ \f1^2\, x^{3} + x^{2} + \f1\, x + 1
$\\\hline
(42, 23, 11)& 92&
$x^{19} + x^{17} + \f1^2\, x^{16} + \f1\, x^{15} + \f1^2\,
x^{14 } + \f1\, x^{13} + \f1\, x^{11} + \f1^2\, x^{10} +
\f1\, x^{9} + x^{7} + \f1\, x^{6} + \f1^2\, x^{5} + \f1\,
x^{4} + \f1\, x + \f1^2
$\\\hline
(40, 16, 15)& 6&
$x^{24} + \f1\, x^{23} + x^{22} + x^{21} + \f1^2\, x^{20} +
\f1\, x^{19} + \f1\, x^{18} + \f1\, x^{17} + x^{15} + x^{14}
+ x^{13} + \f1\, x^{11} + \f1^2\, x^{10} + x^{9} + x^{8} + x^{7}
+ \f1^2\, x^{6} + \f1\, x^{5} + \f1^2\, x^{4} + \f1\,
x^{2} + \f1^2
$\\\hline
(36, 20, 10)& 13&
$x^{16} + \f1^2\, x^{15} + x^{13} + \f1^2\, x^{12} + x^{11} +
\f1\, x^{10} + x^{9} + \f1^2\, x^{8} + \f1\, x^{7} + \f1\,
x^{6} + \f1\, x^{4} + \f1^2\, x^{3} + \f1^2\, x^{2} + 1
$\\\hline
(30,16, 9)& 422& 
$x^{14} + x^{13} + \f1\, x^{11} + x^{10} + x^{9} + x^{8} +
\f1\, x^{7} +x^{6} + \f1\, x^{5} + \f1^2\, x^{4} + \f1^2\,
x^{2} + \f1\, x + \f1^2
$
\end{array}
$$
\caption{Parameters and generating polynomial of skew-cyclic codes over 
$\FF_{4}$.  Here $\alpha$ a zero of $y^2 + y + 1$  and $\theta$ the Frobenius automorphism.
For each code the minimum distance has been improved by 1 according to
Brouwer's table}
\label{codes4}
\end{table}     

For codes over $\FF_{9}$ we managed to improve the lower bound for the 
best known codes in one case (cf. Table~\ref{codes9}).
Due to the larger ground field and the larger codes, the calculation of $d_{{min}}$ is even more 
time consuming than in the previous case. Therefore we stopped our search
for good codes at $n=44$.

\begin{table}
\def\f1{\alpha}
$$\begin{array}{r|r|p{10cm}}
    (n,k,d_{min}) & No & \multicolumn{1}{c}{g}\\ \hline\hline
(44, 20,17)& 5 &
$x^{24} + x^{21} + x^{20} + \f1^7\, x^{19} + \f1^3\, x^{18} +
2\, x^{17} + \f1^3\, x^{16} + \f1^5\, x^{14} + \f1^5\,
x^{13} + 2\, x^{12} + \f1^2\, x^{10} + \f1^7\, x^{9} + 2\,
x^{6} + \f1^5\, x^{5} + \f1^7\, x^{4} + \f1^3\, x^{3} +
\f1^7\, x^{2} + \f1^2\, x + 2
$
\end{array}
$$
\caption{Parameters and generating polynomial of skew-cyclic codes over 
$\FF_{9}$.  Here $\alpha$ a zero of $y^2 - y - 1$ and $\theta$ the Frobenius automorphism.
The minimum distance has been improved by 1 according to the best known 
codes.}
\label{codes9}
\end{table}     

\newpage 
\section{Decoding}

In the following, instead of  a general decoding procedure, we will adapt  (to {\em skew   BCH codes}) the algorithm for decoding BCH codes with designed distance  (see \cite{NS} pages $27-33$ or \cite{WS}) to  $\theta$-cyclic codes.
 We denote $\alpha\in \FF_q$
a primitive $(q-1)$-th root of unity. We suppose that $n$ is even, $q=2^m$ where $m=n$ and  that $\theta(\alpha)=\alpha^2$. Consider  a $\theta$-cyclic code ${\cal C}$ whose generating polynomial is  $G \in \FF_q[X,\theta]$ which is a right divisor of $X^n-1$ in $\FF_q[X,\theta]$. We suppose that  ${\cal C}$
is a skew   BCH codes of designed distance $d\in \NN$, which in this context just means that  $X-\alpha^k$ is a right factor of $G$ for $k\in \{1, \ldots, d-1\}$.  In the following section we give an example of a skew   BCH codes which is not even cyclic in the classical sense, showing that this class extends the class of BCH codes.
The following result allows us to switch to commutative rings for some considerations~:
\begin{pr}
\label{pr1}
For  ${\displaystyle P=\sum_{k=0}^{n-1} a_k \, X^k \in \FF_q[X,\theta]}$, $\beta \in \FF_q$ and $r \in  \FF_q$  the remainder of the right division of $P$ by $X-\beta$, then  $r= \tilde{P}(\beta)$
where $\tilde{P}$ is a (classical) polynomial given by ${\displaystyle \tilde{P}=\sum_{k=0}^{n-1} a_k \, z^{2^k-1} \in  \FF_q[z]}$
\end{pr}
\begin{proof}
The remainder of the right division of $P(X)$ by $X-\beta$ is $$r=a_0 +
a_1 \, \beta + a_2 \, \beta \, \theta(\beta) + a_3 \, \beta \, \theta(\beta) \,
\theta^2(\beta) + \cdots + a_{n-1} \, \beta \, \cdots \theta^{n-2}(\beta)$$
Replacing $\theta^k(\beta)$ with $\beta^{2^k}$, we get ${\displaystyle r=\sum_{k=0}^{n-1}
a_k \, \beta^{2^k-1}=\tilde{P}(\beta)}$
\end{proof}
Therefore the remainder of the right division of $P\in \FF_q[X,\theta]$ by $X-\beta$ (and the image of the remainder in  $\FF_q[X,\theta]/(X^n-1)$) can be interpreted as the evaluation of the polynomial $\tilde{P}$ in the {\em commutative
ring} $\FF_q[z]$ at $\beta\in \FF_q$.\\
Using this property, we can prove like in the classical case (\cite{WS}, Theorem 6.2) that the distance of the code is at least equal to the designed distance $d$.
\begin{pr}
Let $n$ even, $q=2^n$, $\alpha$ a primitive $(q-1)$-th root of unity. Let ${\cal C}$ be a $\theta$-cyclic code with  $\theta(\alpha)=\alpha^2$. Let $G \in \FF_q[X,\theta]$ be its generating polynomial such that $G$ is a right divisor of $X^n-1$ in $\FF_q[X,\theta]$ and  $X-\alpha^k$ is a right factor of $G$ for $k\in \{1, \ldots, d-1\}$.\\
The distance of the code ${\cal C}$ is equal to its designed distance $d$.
\end{pr}
\begin{proof}
According to property (\ref{pr1}), a test matrix for the code is $$H=\begin{matrice}{c}
H_1 \\
H_2
\end{matrice}$$
where 
$$H_1=\begin{matrice}{cccccc}
\alpha_0 & \alpha_1 & \cdots & \alpha_{d-1} & \cdots &\alpha_{n-1}\\
\alpha_0^2 & \alpha_1^2  & \cdots & \alpha_{d-1}^2 & \cdots &\alpha_{n-1}^2\\
\vdots  & & & & & \vdots\\
\alpha_0^{d-1} & & \cdots & & & \alpha_{n-1}^{d-1}
\end{matrice}$$
and $\alpha_i=\alpha^{2^i-1}$. 
If we consider all the possible sets of $d-1$ columns extracted from the $n$ columns of $H_1$, we get square matrices of order $d-1$. Their determinants are non zero if and only if $\alpha_i-\alpha_j$ is non zero for $j<i<n$.\\
But  $\alpha_i-\alpha_j=0 \Leftrightarrow \alpha^{2^i-2^j}=1$ and $0 < 2^i-2^j < 2^n-1$, so as $2^n-1$ is the order of $\alpha$, we get non zero determinants. So each set of $d-1$ columns of $H_1$ are linearly independant, one cannot find any word of weight less than $d$ and the minimum distance of the code is at least $d$.
\end{proof}
We can now adapt {\em almost} entirely the classical decoding algorithm for BCH codes which in described in \cite{NS} pages $27-33$.\\
Let $a\in  \FF_q[X,\theta]/(X^n-1)$ be a code word and let $b=a+e\in \FF_q[X,\theta]/(X^n-1)$ be the received word where
$e=e_{i_1} \, X^{i_1} + \cdots + e_{i_r} \, X^{i_r}$
is the error polynomial with $i_1 < i_2 < \cdots < i_r$ and where $r
\leq t := \frac{d-1}{2}$.\\
One defines the {\em syndrome polynomial} of $e$ as the polynomial
$$S_d(z)=\sum_{k=1}^{d-1} \mbox{Rem}(e,X-\alpha^k) z^{k-1}\in  \FF_q[z].$$
Here the remainder $\mbox{Rem}(e,X-\alpha^k)$ is to be computed in $ \FF_q[X,\theta]$. As $\mbox{Rem}(e,X-\alpha^k)=\mbox{Rem}(b,X-\alpha^k)$, one can compute $S_d(z)$
using the received polynomial $b$.  The syndrome polynomial can
also be written $$S_d(z)=\sum_{k=1}^{d-1} \tilde{e}(\alpha^k) \,
z^{k-1}$$
where $\tilde{e}(z)={\displaystyle \sum_{k=1}^{r} e_{i_k} \, z^{j_k} \in \FF_q[z]}$
 and $j_k=2^{i_k}-1.$\\
One also defines the {\em pseudo-locator polynomial}
$$\sigma(z)=\prod_{k=1}^{r} (1 - \alpha^{j_k} \, z)$$ and the {\em
evaluator polynomial} $$w(z)=\sum_{l=1}^{r} e_{i_l} \, \alpha^{j_l} \,
\prod_{k \neq l} (1-\alpha^{j_k} \, z).$$
Knowing $\sigma(z)$ enables us to find the $j_k$, so that we have {\em almost}
located the positions $i_k$ of the errors in $e$. This point is in fact the only difference with the classical algorithm.

Once we know the
$j_k$ and the evaluator polynomial $w(z)$, we can recover all the $e_{i_k}$ using the following equality
$$e_{i_k}=\alpha^{-j_k} \, w(\alpha^{-j_k}) \, \prod_{l \neq k}
(1-\alpha^{j_l-j_k}), \, k \in \{1, \ldots, r\}.$$
Let us now define $$S(z)=\sum_{k=1}^{\infty} \tilde{e} (\alpha^k) \,
z^{k-1}=S_d(z)+z^{d-1} \, \sum_{k=0}^{\infty}
\tilde{e}(\alpha^{k+1+d}) \, z^{k}$$
Like in \cite{NS} (theorem $I.8$ page $28$), one gets the classical 'key equation':
 $$\sigma(z) \, S(z)=w(z)$$
which one can write
$$\sigma(z) \, S_d(z)+ v(z) \, z^{d-1} = w(z)$$
where ${\displaystyle v(z)=\sigma(z) \, \sum_{k=0}^{\infty} \tilde{e}(\alpha^{k+1+d})
\, z^{k}}$.\\
Following \cite{NS} (theorem $I.11$ page $32$) we apply Euclid's algorithm to the polynomials $S_d(z)$ and
$z^{d-1}$ in $\FF_q[z]$.  We construct the sequences $(r_i(z))$, $(U_i(z))$ and $(V_i(z))$ defined by
$$r_{-1}(z)=z^{d-1}, \; r_0(z)=S_d(z)$$
$$U_{-1}(z)=0, U_0(z)=1, V_{-1}(z)=1, V_0(z)=0$$ and at each step $i$,
$$r_i(z)=r_{i-2}(z)-q_i(z) \, r_{i-1}(z) \; \mbox{with} \; \mbox{deg}(r_{i}(z))< \mbox{deg}(r_{i-1}(z))$$
$$U_i(z)=U_{i-2}(z) - q_i(z) \, U_{i-1}(z), V_i(z)=V_{i-2}(z) - q_i(z)
\, V_{i-1}(z).$$
and we stop as soon as we find $k$
such that $\mbox{deg}(r_{k-1}) \geq t$ and $\mbox{deg}(r_k) < t$.\\
We get
$$U_k(z) \, S_d(z) + V_k(z) \, z^{d-1} = r_k(z),$$
$$\sigma(z)=\frac{U_k(z)}{U_k(0)}$$
and
$$w(z)=\frac{r_k(z)}{r_k(0)}$$
Now from the roots of the pseudo-locator polynomial $\sigma(z)$
we get $j_l, \, l \in \{1, \ldots, r\}$ and from the evaluator
polynomial $w(z)$ we get $$e_{i_l}=\alpha^{-j_l} \, w(\alpha^{-j_l})
\, \prod_{k \neq l} (1-\alpha^{j_k-j_l}), \, l \in \{1, \ldots, r\}.$$
So we have found the coefficients of $e$ and we have almost found the
positions $i_l$ of the errors.  For each $j_l$, we get a finite number
of possibilities $i_l$ solutions to the equation $$j_l \equiv
2^{i_l}-1 \pmod{n}$$
So we get a finite number of possible errors, which we test until we
find $e$ such that $b+e$ is a code word. As the distance of the code is $d$ we are sure that such a $e$ is unique and so we have decoded.\\

\section{Worked example}

Let $n=m=10$ and let $\alpha$ such that $\alpha^{2^{10}-1}=1$.  The
polynomial $$G=X^6 + \alpha^{345} \, X^5 + \alpha^{643} \, X^4 +
\alpha^{878} \, X^3 + \alpha^{670} \, X^2 + \alpha^{1020} \, X +
\alpha^{777}$$
divides $X^{10}+1$ to the right in $\FF_{2^{10}}[X,\theta]$.  Therefore it is the
generator polynomial of a $\theta$-cyclic code ${\cal C}$ of length $10$ over
$\FF_{2^{10}}$. Since $X-\alpha^k$ is a right factor of $G$ for $k\in \{1, \ldots, 6\}$, the code 
${\cal C}$ is of  designed distance $d=7$. One can check that this skew BCH code is not cyclic in the classical sense.

We consider the  code word $a$ given by
\begin{center}
  {\small $a(X)=\alpha^{654}\,X^9 + \alpha^{547}\,X^8 +
\alpha^{650}\,X^7 + \alpha^{16}\,X^6 + \alpha^{567}\,X^5 +
\alpha^{29}\,X^4 + \alpha^{87}\,X^3 + \alpha^{696}\,X^2 +
\alpha^{252}\,X + \alpha^{555},$}
\end{center}
an error 
$$e=\alpha^{341}\,X^9 + \alpha^{682}\,X^8 +
\alpha^{682}.$$ 
The received perturbed code word $b=a+e$  is 
\begin{center}
{\small $b=\alpha^{818}X^9 +
\alpha^{775}X^8 + \alpha^{650}X^7 + \alpha^{16}X^6 +
\alpha^{567}X^5 + \alpha^{29}X^4 + \alpha^{87}X^3 +
\alpha^{696}X^2 + \alpha^{252}X + \alpha^{557}.$}
\end{center}
Knowing the received polynomial $b$ and $d=7$, we can compute the
syndrome polynomial $$S_7(z)=\alpha^{404}\,z^5 + \alpha^{403}\,z^4 +
\alpha^{601}\,z^3 + \alpha^{645}\,z^2 + \alpha^{614}\,z +
\alpha^{406}$$
Applying Euclid algorithm to $S_7(z)$ and $z^6$ in $\FF_{2^{10}}[z]$
with $t=3$, we get the {\em pseudo-locator polynomial} $\sigma(z)$
$$\sigma(z)=\alpha^{766}\,z^3 + \alpha^{642}\,z^2 + \alpha^{241}\,z +
1$$
and the {\em evaluator polynomial} $w(z)$ $$w(z)=\alpha^{84}\,z^2 +
\alpha^{185}\,z + \alpha^{406}.$$
From the roots $1, \alpha^{512}$ and $\alpha^{768}$ of the polynomial
$\sigma(z)$ we get the value of $r$ ($r=3$) and the values of $j_1,
j_2, j_3$~: $$j_1=0, j_2=511, j_3=255.$$
We can now find the values of the coefficients of $e$ via the
polynomial $w$~:
$$
    e_{i_1}=\alpha^{682}, e_{i_2}= \alpha^{341}, e_{i_3}=
    \alpha^{682}.
$$
We have now to locate exactly the positions of the errors.  For each
$k$ in $\{1,2,3\}$, we solve the equations $$2^{i_k} - 1 \equiv j_k
\pmod{10}.$$
$$\begin{array}{lcl} 2^{i_1} - 1 \equiv 0 \pmod{10}& \Leftrightarrow & i_1
\, \pmod{10} = 0\\
2^{i_2} - 1 \equiv 511 \pmod{10}&  \Leftrightarrow & i_2 \, \pmod{10} \in
\{1,5,9\}\\
2^{i_3} - 1 \equiv 255 \pmod{10} & \Leftrightarrow  & i_3 \, \pmod{10} \in
\{4,8\} \end{array}$$
So the list of the possible errors is $$\begin{array}{ll}
[\alpha^{682}\,X^4 + \alpha^{341}\,X + \alpha^{682}, &
    \alpha^{682}\,X^8 + \alpha^{341}\,X + \alpha^{682},\\
    \alpha^{341}\,X^5 + \alpha^{682}\,X^4 + \alpha^{682},&
    \alpha^{682}\,X^8 + \alpha^{341}\,X^5 + \alpha^{682},\\
    \alpha^{341}\,X^9 + \alpha^{682}\,X^4 + \alpha^{682},&
    \alpha^{341}\,X^9 + \alpha^{682}\,X^8 + \alpha^{682}]
    \end{array}$$
The only one such that $g$ divides $b+e$ to the right is
$$e=\alpha^{341}\,X^9 + \alpha^{682}\,X^8 + \alpha^{682}.$$

For this code, $5000$ random tests have been made (each random test
takes a random code word, a random error of weight at most three and
checks whether the corrected word is equal to the code word).\\

\end{document}